\documentclass{article}
\usepackage{color}
\catcode`\@=11 \@addtoreset{equation}{section}

\catcode`\@=12

\newtheorem{Theorem}{Theorem}[section]
\newtheorem{Proposition}{Proposition}[section]
\newtheorem{Lemma}{Lemma}[section]
\newtheorem{Corollary}{Corollary}[section]

\newcommand{\bTheorem}[1]{
\begin{Theorem} \label{T#1} }
\newcommand{\eT}{\end{Theorem}}

\newcommand{\bProposition}[1]{
\begin{Proposition} \label{P#1}}
\newcommand{\eP}{\end{Proposition}}

\newcommand{\bLemma}[1]{
\begin{Lemma} \label{L#1} }
\newcommand{\eL}{\end{Lemma}}

\newcommand{\bCorollary}[1]{
\begin{Corollary} \label{C#1} }
\newcommand{\eC}{\end{Corollary}}

\newcommand{\bFormula}[1]{
\begin{equation} \label{#1}}
\newcommand{\eF}{\end{equation}}

\newcommand{\Ov}[1]{\overline{#1}}
\newcommand{\DC}{C^\infty_c}

\newcommand{\vr}{\varrho}
\newcommand{\vre}{\vr_\ep}
\newcommand{\vte}{\vt_\ep}
\newcommand{\vue}{\vu_\ep}

\newcommand{\vt}{\vartheta}
\newcommand{\vu}{\vc{u}}
\newcommand{\vc}[1]{{\bf #1}}

\newcommand{\Div}{{\rm div}_x}
\newcommand{\Grad}{\nabla_x}

\newcommand{\tn}[1]{\mbox {\F #1}}
\newcommand{\dx}{{\rm d} {x}}
\newcommand{\dt}{{\rm d} t }

\newcommand{\dxdt}{\dx \ \dt}
\newcommand{\intO}[1]{\int_{R^3} #1 \ \dx}

\newcommand{\ep}{\varepsilon}

\font\F=msbm10 scaled 1000

\definecolor{grey}{rgb}{0.85,0.85,0.85}
\date{}

\long\def\greybox#1{%
    \newbox\contentbox%
    \newbox\bkgdbox%
    \setbox\contentbox\hbox to \hsize{%
        \vtop{
            \kern\columnsep
            \hbox to \hsize{%
                \kern\columnsep%
                \advance\hsize by -2\columnsep%
                \setlength{\textwidth}{\hsize}%
                \vbox{
                    \parskip=\baselineskip
                    \parindent=0bp
                    #1
                }%
                \kern\columnsep%
            }%
            \kern\columnsep%
        }%
    }%
    \setbox\bkgdbox\vbox{
        \color{grey}
        \hrule width  \wd\contentbox %
               height \ht\contentbox %
               depth  \dp\contentbox
        \color{black}
    }%
    \wd\bkgdbox=0bp%
    \vbox{\hbox to \hsize{\box\bkgdbox\box\contentbox}}%
    \vskip\baselineskip%
}

\begin{document}


\title{Inviscid incompressible limits of the full Navier-Stokes-Fourier system}
\author{Eduard Feireisl \thanks{The work was supported by Grant 201/09/
0917 of GA \v CR and by  RVO: 67985840} \and Anton\' \i n Novotn\' y \thanks{The work was partially supported by  RVO: 67985840}}

\maketitle

\bigskip

\centerline{Institute of Mathematics of the Academy of Sciences of
the Czech Republic,} \centerline{\v Zitn\' a 25, 115 67 Praha 1,
Czech Republic} \centerline{and} \centerline{Erwin Schroedinger
International Institute for Mathematical Physics,} \centerline{
Boltzmanngasse 9, A-1090 Vienna, Austria}

\bigskip

\centerline{IMATH, Universit\' e du Sud Toulon-Var}

\centerline{\v BP 20139, 839 57 La Garde, France}

\medskip

\begin{abstract}
We consider the full Navier-Stokes-Fourier system in the singular limit for the small
Mach and large Reynolds and P{\' e}clet numbers, with ill prepared initial data on $R^3$.
The Euler-Boussinesq approximation is identified as the limit system.
\end{abstract}

\bigskip


\medskip


\bigskip

\tableofcontents

\section{Introduction}
\label{i}

Scale analysis and the associated mathematical problems of singular limits reveal the dominant features of complete fluid systems in the regime where some characterictic parameters become small or infinitely large. We apply the method
of \emph{relative entropies} developed in \cite{FeiNov10} to study the asymptotic limit in the complete Navier-Stokes-Fourier system for low Mach and large Reynolds and P{\'e}clet numbers.

In order to avoid unnecessary technical difficulties, we consider the hypothetical situation when a compressible fluid, described by means of the Eulerian density $\vr = \vr(t,x)$, the velocity field $\vu = \vu(t,x)$, and the absolute temperature $\vt = \vt(t,x)$ occupies the entire physical space $R^3$. The associated Navier-Stokes-Fourier system of field equations reads:

\bFormula{i1}
\partial_t \vr + \Div (\vr \vu) = 0,
\eF
\bFormula{i2}
\partial_t (\vr \vu) + \Div (\vr \vu \otimes \vu) + \frac{1}{\ep^2} \Grad p(\vr, \vt) = \ep^a \Div \tn{S} (\vt, \Grad \vu),
\eF
\bFormula{i3}
\partial_t (\vr s(\vr , \vt)) + \Div (\vr s(\vr, \vt) \vu) + \ep^\beta \Div \left( \frac{\vc{q}(\vt, \Grad \vt)}{\vt} \right) =
\frac{1}{\vt} \left( \ep^{2 + a} \tn{S} (\vt, \Grad \vu) : \Grad \vu - \ep^b \frac{\vc{q}(\vt, \Grad \vt) \cdot
\Grad \vt }{\vt} \right),
\eF
where $p = p(\vr,\vt)$ is the pressure, $s = s(\vr, \vt)$ the specific entropy, while the symbol $\tn{S}(\vt, \Grad \vu)$ denotes the viscous stress satisfying \emph{Newton's rheological law}
\bFormula{i3a}
\tn{S}(\vt, \Grad \vu) = \mu(\vt) \left( \Grad \vu + \Grad^t \vu - \frac{2}{3} \Div \vu \right),
\eF
and $\vc{q} = \vc{q}(\vt, \Grad \vt)$ is the heat flux determined by \emph{Fourier's law}
\bFormula{i3b}
\vc{q}(\vt, \Grad \vt) = - \kappa (\vt) \Grad \vt.
\eF
Note that, again for the sake of simplicity and clarity of presentation, we have omitted the effect of any external force in the momentum equation (\ref{i2}) as well as the bulk viscosity contribution to the viscous stress (\ref{i3a}).

The scaling of the pressure in (\ref{i2}) corresponds to the Mach number proportional to a small parameter $\ep$, whereas the Reynolds and P{\' e}clet numbers scale as $\ep^{-a}$ and $\ep^{-b}$, respectively. We consider the
initial data in the form
\bFormula{i4}
\vr(0, \cdot) = \vr_{0,\ep} = \Ov{\vr} + \ep \vr^{(1)}_{0,\ep},\
\vt(0, \cdot) = \vt_{0,\ep} = \Ov{\vt} + \ep \vt^{(1)}_{0,\ep},\ \vu(0, \cdot) = \vu_{0, \ep},
\eF
together with the boundary conditions ``at infinity''
\bFormula{i5}
\vr \to \Ov{\vr}, \ \vt \to \Ov{\vt}, \ \vu \to 0 \ \mbox{as} \ |x| \to \infty.
\eF

Under these circumstances, the limit (target) problem can be identified as the incompressible Euler system
\bFormula{l1}
\Div \vc{v} = 0,
\eF
\bFormula{l2}
\partial_t \vc{v} + \vc{v} \cdot \Grad \vc{v} + \Grad \Pi = 0,
\eF
supplemented with a transport equation for the temperature deviation $T$,
\bFormula{l3}
\partial_t T + \vc{v} \cdot \Grad T = 0.
\eF
Here, the function $\vc{v}$ is the limit velocity while $T \approx \frac{\vt - \Ov{\vt}}{\ep}$. Note that the system (\ref{l1} - \ref{l3}) can be obtained
as a hydrodynamic of the Boltzmann equation, see Golse \cite{Golse05}.
The exact statement of our results including the initial data for the target system (\ref{l1} - \ref{l3}) will be given in Theorem \ref{Tm1} below.

Our approach is based on the concept of (very) weak solutions for the Navier-Stokes-Fourier system (\ref{i1} - \ref{i3}), developed in \cite{FeNo6}, and extended to problems on unbounded domains in \cite{JeJiNo}. Accordingly, the convergence to the limit problem takes place on any time interval $[0,T]$ on which the Euler system (\ref{l1}), (\ref{l2}) possesses a regular solution. Similar results for the compressible barotropic Navier-Stokes system were obtained by Masmoudi \cite{MAS5}, see also the survey paper \cite{MAS1} of the same author.
Alazard \cite{AL2}, \cite{AL1}, \cite{AL} studies the singular limits of the compressible Euler and the Navier-Stokes-Fourier system using the approach proposed by Klainerman and Majda \cite{KM1} based on {strong solutions}.  To the best of our knowledge, the present paper represents the first result of this kind for compressible \emph{and} heat conducting fluids in the framework of weak solutions. The main novelty of our approach is the use of the relative entropy for the Navier-Stokes-Fourier system
discovered in \cite{FeiNov10} to establish the necessary \emph{uniform} bounds independent of the scaling parameter $\ep$, and, more importantly, to obtain \emph{stability} of solutions to the limit system.

The paper is organized as follows. In Section \ref{v}, we collect the necessary preliminary material and introduce the concept of \emph{very weak solution} to the Navier-Stokes-Fourier system on $R^3$. The main result on the asymptotic limit for $\ep \to 0$ is stated in Section \ref{m}. The remaining part of the paper will be devoted to the proof of the main theorem. In Section \ref{b}, we use the total dissipation balance associated to the Navier-Stokes-Fourier system to establish all necessary uniform bounds independent of $\ep \to 0$. The crucial ingredient of the proof is the
relative entropy inequality introduced in Section \ref{r} that provides the necessary stability estimates for the limit system. As is typical for this kind of problems, the most difficult part is to establish the convergence of the oscillatory gradient component of the velocity field
corresponding to the presence of acoustic waves. Since the problem is considered on the whole space $R^3$, this can be accomplished by the standard dispersive estimates, see Section \ref{a}. Finally, the proof of convergence towards the limit system is finished in Section \ref{c}.

\section{Preliminaries, very weak solutions for the full Navier-Stokes-Fourier system}
\label{v}

We start by listing the technical hypotheses imposed on constitutive relations. They are analogous to those introduced
in the framework of the existence theory developed in \cite[Chapter 3]{FeNo6}, where the interested reader can find all relevant information concerning the physical background as well as possible generalizations.

We suppose that the pressure $p=p(\vr,\vt)$ is given by the formula
\bFormula{i10}
p(\vr, \vt) = \vt^{5/2} P \left( \frac{\vr}{\vt^{3/2}} \right) + \frac{a}{3} \vt^4 , \ a  > 0,
\eF
while the specific internal energy $e = e(\vr,\vt)$ and the specific entropy $s = s(\vr, \vt)$ read
\bFormula{i11}
e(\vr, \vt) = \frac{3}{2} \vt \frac{ \vt^{3/2} }{\vr} P \left( \frac{\vr}{\vt^{3/2}} \right) + {a} \vt^4
\eF
\bFormula{i12}
s(\vr, \vt) = S \left( \frac{\vr}{\vt^{3/2}} \right) + \frac{4a}{3} \frac{\vt^3}{\vr},
\eF
where
\bFormula{i13}
P \in C^1[0, \infty) \cap C^3(0,\infty), \ P(0) = 0 ,\ P'(Z) > 0 \ \mbox{for all}\ Z \geq 0,
\eF
\bFormula{i14}
\lim_{Z \to \infty} \frac{P(Z)}{Z^{5/3}} = P_\infty > 0,
\eF
\bFormula{i15}
0 < \frac{ \frac{5}{3} P(Z) - P'(Z) Z }{Z} < c \ \mbox{for all} \ Z > 0,
\eF
and
\bFormula{i16}
S'(Z) = - \frac{3}{2} \frac{ \frac{5}{3} P(Z) - P'(Z) Z }{Z^2} , \ \lim_{Z \to \infty} S(Z) = 0.
\eF
Let us only remark that the rather misteriously looking relation (\ref{i15}) expresses positivity and uniform boundedness of the specific heat at constant volume.

In addition, the transport coefficients $\mu$ and $\kappa$ vary with the temperature, specifically,
\bFormula{i8}
\mu \in C^1[0, \infty) \ \mbox{is (globally) Lipschitz continuous,}\ 0 < \underline{\mu} (1 + \vt) \leq \mu(\vt)
\ \mbox{for all}\ \vt \geq 0,
\eF
\bFormula{i9}
\kappa \in C^1[0, \infty),\ 0 < \underline{\kappa}(1 + \vt^3) \leq \kappa (\vt) \leq \Ov{\kappa}(1 + \vt^3) \ \mbox{for all}\ \vt \geq 0.
\eF

\subsection{Very weak solutions}

To begin, we introduce the \emph{ballistic free energy}
\bFormula{bfe}
H_\Theta (\vr, \vt) = \vr \Big( e(\vr, \vt) - \Theta s(\vr,\vt) \Big), \ \mbox{where}\  \Theta > 0.
\eF
Following \cite{JeJiNo} we shall say that a trio of functions $\{ \vr, \vt, \vu \}$ represents a \emph{very weak solution} of the Navier-Stokes-Fourier system (\ref{i1} - \ref{i5}) on the space time cylinder $(0,T) \times \Omega$ if:

\begin{itemize}

\item $\vr \geq 0$, $\vt > 0$ a.a. in $(0,T) \times \Omega$,
\[
( \vr - \Ov{\vr} ) \in L^\infty(0,T; L^2 + L^{5/3}(R^3),\ (\vt - \Ov{\vt}) \in L^\infty(0, T; L^2 + L^4 (R^3)),
\]
\[
\Grad \vt , \ \Grad \log(\vt) \in L^2(0,T; L^2(R^3;R^3)),
\]
\[
\vu \in L^2(0,T; W^{1,2}(R^3; R^3));
\]

\item the equation of continuity (\ref{i1}) is replaced by a family of integral identities
\bFormula{v1}
\int_{R^3} \Big[ \vr(\tau, \cdot) \varphi (\tau, \cdot) - \vr_{0,\ep} \varphi (0, \cdot) \Big] \ \dx
= \int_0^\tau \int_{R^3} \Big( \vr \partial_t \varphi + \vr \vu \cdot \Grad \varphi \Big) \ \dxdt
\eF
for any $\tau \in [0,T]$ and any test function $\varphi \in \DC([0,T] \times R^3)$;

\item the momentum equation (\ref{i2}), together with the initial condition (\ref{i4}), is satisfied in the sense of distributions, specifically,
\bFormula{v2}
\int_{R^3} \Big[ \vr \vu (\tau, \cdot) \cdot \varphi (\tau, \cdot) - \vr_{0,\ep} \vu_{0,\ep}
\varphi(0, \cdot) \Big] \ \dx
\eF
\[
= \int_0^\tau \int_{R^3} \Big( \vr \vu \cdot \partial_t \varphi + \vr \vu \otimes \vu : \Grad \varphi +
\frac{1}{\ep^2} p(\vr,\vt) \Div \varphi - \ep^{a} \tn{S}(\vt, \Grad \vu) : \Grad \varphi \Big) \ \dxdt
\]
for any $\tau \in [0,T]$, and any $\varphi \in \DC([0,T] \times R^3; R^3)$;

\item the entropy production equation (\ref{i3}) is relaxed to the entropy inequality
\bFormula{v3}
\int_{R^3} \Big[ \vr_{0,\ep} s(\vr_{0,\ep}, \vt_{0, \ep} ) \varphi(0, \cdot) -
\vr s(\vr, \vt) (\tau, \cdot) \varphi(\tau, \cdot) \Big] \ \dx
\eF
\[
+ \int_0^\tau \int_{R^3} \frac{1}{\vt} \left( \ep^{2 + a} \tn{S}(\vt, \Grad \vu) : \Grad \vu - \ep^b
\frac{\vc{q}(\vt, \Grad \vt) \cdot \Grad \vt }{\vt} \right) \varphi \ \dxdt
\]
\[
\leq - \int_0^\tau \int_{R^3} \left( \vr s(\vr,\vt) \partial_t \varphi + \vr s(\vr, \vt) \vu \cdot \Grad \varphi +
\ep^b \frac{ \vc{q}(\vt, \Grad \vt) }{\vt} \cdot \Grad \varphi \right) \ \dxdt
\]
for a.a. $\tau \in [0,T]$ and any test function $\varphi \in \DC([0,T] \times R^3)$, $\varphi \geq 0$;

\item the \emph{total dissipation inequality}
\bFormula{v4}
\int_{R^3} \left[ \frac{1}{2} \vr |\vu|^2 + \frac{1}{\ep^2} \left( H_{\Ov{\vt}}(\vr, \vt) -
\frac{\partial H_{\Ov{\vt}}(\Ov{\vr}, \Ov{\vt}) }{\partial \vr}(\vr - \Ov{\vr}) - H_{\Ov{\vt}} (\Ov{\vr}, \Ov{\vt})
\right) (\tau, \cdot) \right] \ \dx
\eF
\[
+ \Ov{\vt} \int_0^\tau \int_{R^3} \frac{1}{\vt} \left( \ep^{a} \tn{S}(\vt, \Grad \vu) : \Grad \vu - \ep^{b - 2} \frac{\vc{q}(\vt, \Grad \vt) \cdot \Grad \vt }{\vt} \right) \ \dxdt
\]
\[
\leq \int_{R^3} \left[ \frac{1}{2} \vr_{0,\ep} |\vu_{0,\ep}|^2 + \frac{1}{\ep^2} \left( H_{\Ov{\vt}}(\vr_{0,\ep},
\vt_{0,\ep}) -
\frac{\partial H_{\Ov{\vt}}(\Ov{\vr}, \Ov{\vt}) }{\partial \vr}(\vr_{0,\ep} - \Ov{\vr}) - H_{\Ov{\vt}} (\Ov{\vr}, \Ov{\vt})
\right)  \right] \ \dx
\]
holds for a.a. $\tau \in [0,T]$.

\end{itemize}

Under the hypotheses (\ref{i10} - \ref{i9}), the existence of very
weak solutions to the Navier-Stokes-Fourier system in $(0,T)
\times R^3$ was shown in \cite{JeJiNo}, {  along with the property
that a very weak solution coincides with the strong solution
emanating from the same initial as long as the latter exists
(known as the weak-strong uniqueness principle).}

\section{Main result}

\label{m}

Suppose that $\vc{v}_0$ is a given vector field such that
\[
\vc{v}_0 \in W^{k,2}(R^3;R^3),\ k > \frac{5}{2}, \ \|\vc{v}_0\|_{W^{k,2}(\Omega;R^3)}\le D, \ \Div \vc{v}_0 = 0.
\]
It is well-known that the Euler system (\ref{l1}), (\ref{l2}), supplemented with the initial condition
\[
\vc{v}(0, \cdot) = \vc{v}_0.
\]
possesses a regular solution $\vc{v}$, unique in the class
\bFormula{euler} \vc{v} \in C([0, T_{\rm max}); W^{k,2}(R^3;R^3)),
\ \partial_t \vc{v} \in C([0, T_{\rm max}); W^{k-1,2}(R^3;R^3)),
\eF defined on a maximal time interval $[0, T_{\rm max})$, $T_{\rm
max} = T_{\rm max}(D)$, see Kato \cite{Kato}.

For each vector field $\vc{U} \in L^2(R^3;R^3)$ we denote by $\vc{H}[\vc{U}]$ the standard \emph{Helmholtz projection}
on the space of solenoidal functions.

We are ready to state the main result of this paper.

\greybox{

\bTheorem{m1}
Let the thermodynamic functions
$p$, $e$, and $s$ comply with hypotheses (\ref{i10} - \ref{i16}), and let the transport coefficients $\mu$ and
$\kappa$  satisfy (\ref{i8}), (\ref{i9}). Let
\bFormula{cof}
b > 0, \ 0 < a < \frac{10}{3}.
\eF
Furthermore, suppose that
the initial data (\ref{i4}) are chosen in such a way that
\bFormula{m1}
\{ \vr^{(1)}_{0, \ep} \}_{\ep > 0},\  \{ \vt^{(1)}_{0,\ep} \}_{\ep > 0} \ \mbox{are bounded in}\
L^2 \cap L^\infty (R^3), \ \vr^{(1)}_{0, \ep} \to \vr^{(1)}_0,\ \vt^{(1)}_{0, \ep} \to \vt^{(1)}_0
\ \mbox{in}\ L^2(R^3),
\eF
and
\bFormula{m2}
\{ \vu_{0,\ep} \}_{\ep > 0} \ \mbox{is bounded in}\ L^2(R^3;R^3), \ \vu_{0, \ep} \to \vu_0 \ \mbox{in}\ L^2(R^3;R^3),
\eF
where
\bFormula{HYP}
\vr^{(1)}_0 , \ \vt^{(1)}_0 \in W^{1,2} \cap W^{1,\infty}(R^3),\ \vc{H}[\vu_0] = \vc{v}_0 \in W^{k,2}(R^3;R^3) \ \mbox{for a certain}\ k > \frac{5}{2}.
\eF
Let $T_{\rm max} \in (0, \infty]$ denote the maximal life-span of the regular solution $\vc{v}$ to the Euler system (\ref{l1}), (\ref{l2}) satisfying $\vc{v}(0, \cdot) = \vc{v}_0$. Finally, let $\{ \vre, \vte, \vue \}$ be a very weak solution of the Navier-Stokes-Fourier system (\ref{i1} - \ref{i5}) in $(0,T) \times R^3$, $T < T_{\rm max}$.

Then
\bFormula{m3}
{\rm ess} \sup_{t \in (0,T)} \| \ \vre (t, \cdot) - \Ov{\vr} \ \|_{L^2 + L^{5/3} (R^3)} \leq \ep c,
\eF
\bFormula{m4}
\sqrt{ {\vre} } \vue \to \sqrt{ \Ov{\vr} } \ \vc{v} \ \mbox{in}\ L^\infty_{\rm loc}((0,T]; L^2_{\rm loc} (R^3; R^3))
\ \mbox{and weakly-(*) in} \ L^\infty(0,T; L^2(R^3;R^3)),
\eF
and
\bFormula{m5}
\frac{ \vte - \Ov{\vt} }{\ep} \to T \ \mbox{in}\ L^\infty_{\rm loc}((0,T]; L^q_{\rm loc} (R^3; R^3)),\ 1 \leq q < 2,
\ \mbox{and weakly-(*) in}\ L^\infty(0,T; L^2(R^3)),
\eF
where $\vc{v}$, $T$ is the unique solution of the Euler-Boussinesq system (\ref{l1} - \ref{l3}), with the initial data
\bFormula{m6}
\vc{v}_0 = \vc{H}[ \vu_0 ], \ T_0 = \Ov{\vr} \frac{ \partial s(\Ov{\vr}, \Ov{\vt}) }{\partial \vt} \vt^{(1)}_{0} -
\frac{1}{\Ov{\vr}} \frac{ \partial p(\Ov{\vr}, \Ov{\vt}) }{\partial \vt} \vr^{(1)}_{0}.
\eF

\eT

}

It is worth noting that the initial distribution of the temperature deviation $T_0$ includes a contribution proportional to
$\vr^{(1)}_{0}$. This is related to the well-known data adjustment problem observed by physicists, see Zeytounian \cite{ZEY2} and the discussion in \cite[Chapter 5.5.3]{FeNo6}.
The rest of the paper is devoted to the proof of Theorem \ref{Tm1}.

\section{Uniform bounds}
\label{b}

Thanks to the hypotheses (\ref{m1}), (\ref{m2}), the integral on the right-hand side of the total dissipation inequality (\ref{v4}) remains bounded uniformly for $\ep \to 0$. On the other hand, in accordance with the structural properties of the thermodynamic functions stated in (\ref{i10} - \ref{i16}), the function
\[
 {H_{\Theta}
(\vr, \vt)} - \frac{ \partial H_{\Theta} (r, \Theta)}{\partial \vr}
(\vr - r) - H_{\Theta} ( r , \Theta )
\]
enjoys the following coercivity properties:
For any compact $K \subset (0, \infty)^2$ and
\[
(r, \Theta) \in K,
\]
there exists a strictly positive constant $c(K)$, depending only on $K$ and the structural properties of $P$, such that
\bFormula{b1}
{H_{\Theta}
(\vr, \vt)} - \frac{ \partial H_{\Theta} (r, \Theta)}{\partial \vr}
(\vr - r) - H_{\Theta} ( r , \Theta ) \geq c(K) \left( |\vr - r |^2 + |\vt - \Theta |^2 \right) \ \mbox{if}
\ (\vr, \vt) \in K,
\eF
\bFormula{b2}
{H_{\Theta}
(\vr, \vt)} - \frac{ \partial H_{\Theta} (r, \Theta)}{\partial \vr}
(\vr - r) - H_{\Theta} ( r , \Theta )
\eF
\[
\geq c(K) \Big( \vr e(\vr,\vt) + \vr \Theta |s(\vr, \vt)| + 1  \Big) \ \mbox{if}
\ (\vr, \vt) \in (0,\infty)^2 \setminus K,
\]
see \cite[Proposition 3.2]{FeNo6}.

In view of (\ref{b1}), (\ref{b2}) it is convenient to introduce a decomposition
\[
h = [h]_{\rm ess} + [h]_{\rm res} \ \mbox{for a measurable function} \ h,
\]
where
\[
[h]_{\rm ess} = h \ 1_{ \{ \Ov{\vr} / 2 < \vre < 2 \Ov{\vr} ; \ \Ov{\vt}/2 < \vte < 2 \Ov{\vt} \} },\
[h]_{\rm res} = h - h_{\rm ess},
\]
see \cite[Chapter 4.7]{FeNo6}.

Consequently, combining (\ref{v4}) with (\ref{b1}), (\ref{b2}) and the hypotheses (\ref{i10} - \ref{i9}) we deduce the following list of estimates:
\bFormula{b3}
{\rm ess} \sup_{t \in (0,T)} \| \sqrt{\vre} \vue (t, \cdot) \|_{L^2(R^3;R^3)} \leq c,
\eF
\bFormula{b4}
{\rm ess} \sup_{t \in (0,T)} \left\| \left[ \frac{\vre - \Ov{\vr} }{\ep}(t, \cdot) \right]_{\rm ess} \right\|_{L^2(R^3;R^3)} + {\rm ess} \sup_{t \in (0,T)} \left\| \left[ \frac{\vte - \Ov{\vt} }{\ep} (t, \cdot) \right]_{\rm ess} \right\|_{L^2(R^3;R^3)} \leq c,
\eF
\bFormula{b5}
{\rm ess} \sup_{t \in (0,T)} \int_{R^3} \left( \left[ \vre^{5/3}(t, \cdot) \right]^{5/3}_{\rm res} + \left[ \vte (t, \cdot) \right]^4_{\rm res} + 1_{\rm res}(t, \cdot) \right) \ \dx \leq \ep^2 c,
\eF
and
\bFormula{b6}
\left\| \ep^{a/2} \vue \right\|_{L^2(0,T; W^{1,2}(R^3;R^3))} \leq c,
\eF
\bFormula{b7}
\left\| \ep^{(b-2)/2} \left( \vte - \Ov{\vt} \right) \right\|_{L^2(0,T; W^{1,2}(R^3;R^3))} + \left\| \ep^{(b-2)/2} \left( \log(\vte) - \log(\Ov{\vt}) \right) \right\|_{L^2(0,T; W^{1,2}(R^3;R^3))} \leq c,
\eF
where the symbol $c$ denotes a generic constant independent of $\ep$. We remark that (\ref{b6}) follows from the generalized Korn's inequality
\[
\left\| \Grad \vc{w} + \Grad^t \vc{w} - \frac{2}{3} \Div \vc{w} \tn{I} \right\|_{L^2(R^3)} \geq c
\left\| \Grad \vc{w}  \right\|_{L^2(R^3)} \ \mbox{for}\ \vc{w} \in W^{1,2}(R^3;R^3),
\]
combined with the estimates (\ref{b3}), (\ref{b5}). Similar arguments based on the Sobolev inequality and (\ref{b4}), (\ref{b5}) yield (\ref{b7}).

\section{Relative entropy inequality}
\label{r}

Motivated by \cite{FeiNov10}, we introduce the
\emph{relative entropy}
\bFormula{r0}
\mathcal{E}_\ep \left( \vr, \vt, \vu \Big| r , \Theta, \vc{U} \right) =
\intO{ \left[ \frac{1}{2} \vr |\vu - \vc{U} |^2 + \frac{1}{\ep^2}
\left( H_\Theta (\vr, \vt) - \frac{\partial H_\Theta (r, \Theta) }{\partial \vr} (\vr - r) -
H_\Theta (r, \Theta) \right) \right] },
\eF
where $H_\Theta$ was defined through (\ref{bfe}). As shown in \cite{JeJiNo}, any very weak solution
solution $\{ \vr, \vt, \vu \}$ of the scaled Navier-Stokes-Fourier system satisfies the \emph{relative entropy inequality} in the form:

\greybox{

\bFormula{r1}
\left[ \mathcal{E}_\ep \left( \vr, \vt, \vu \Big| r , \Theta, \vc{U} \right) \right]_{t = 0}^\tau
+ \int_0^\tau \intO{ \frac{\Theta}{\vt} \left( \ep^a \tn{S} (\vt, \Grad \vu) : \Grad \vu -
\ep^{b-2} \frac{\vc{q}(\vt, \Grad \vt) \cdot \Grad \vt }{\vt} \right) } \ \dt
\eF
\[
\leq \int_0^\tau \intO{ \Big( \vr \Big( \partial_t \vc{U} +
\vu\cdot \Grad \vc{U} \Big) \cdot (\vc{U} - \vu)  + \ep^a
\tn{S}(\vt, \Grad \vu): \Grad \vc{U} \Big) } \ \dt
\]
\[
+\frac
1{\ep^2}\int_0^\tau\intO{\Big[\Big(p(r,\Theta)-p(\vr,\vt)\Big){\rm
div}\vc U +\frac\vr {r}(\vc U-\vu)\cdot\Grad
p(r,\Theta)\Big]}{\rm d}t
\]
\[
- \frac{1}{\ep^2} \int_0^\tau \intO{ \left( \vr \Big( s(\vr,\vt) - s(r, \Theta) \Big) \partial_t \Theta +
\vr \Big( s(\vr,\vt) - s(r, \Theta) \Big) \vu \cdot \Grad \Theta + \ep^b \frac{\vc{q}(\vt, \Grad \vt) }{\vt} \cdot
\Grad \Theta \right) } \ \dt
\]
\[
+ \frac{1}{\ep^2} \int_0^\tau \intO{ \frac{r - \vr}
{r}\Big(
\partial_t p(r, \Theta) + \vc U \cdot \Grad p(r,
\Theta)\Big)  } \ \dt
\]

}

\noindent
for any trio of continuously differentiable ``test'' functions defined on $[0,T] \times R^3$,
\[
r > 0 , \ \Theta > 0 , \ r \equiv \Ov{\vr} , \ \Theta \equiv \Ov{\vt} \ \mbox{outside a compact subset of}\ R^3,
\]
\[
\vc{U} \in C([ 0, T]; W^{k,2}(R^3;R^3)),\ \partial_t \vc{U} \in \ C([0,T]; W^{k-1,2}(R^3;R^3)),\ k > \frac{5}{2}.
\]
It seems interesting to notice that the mere relative entropy inequality (\ref{r1}) could be taken as a definition of ``dissipative'' solutions to the Navier-Stokes-Fourier system in the spirit of a similar concept introduced by Lions \cite[Chapter 4.4]{LI4a} in the context of the incompressible Euler system.

We take
\[
\vr = \vre,\ \vt = \vte, \ \vu = \vue
\]
and choose the functions $\{ r, \Theta, \vc{U} \}$ in the
following way:
\bFormula{A}
r = r_\ep = \Ov{\vr} + \ep R_\ep ,\
\Theta = \Theta_\ep = \Ov{\vt} + \ep T_\ep , \ \vc{U} = \vc{U}_\ep
= \vc{v} + \Grad \Phi_\ep; \eF
where $\vc{v}$ is the solution to
the incompressible Euler system (\ref{l1}), (\ref{l2}), with the
initial condition (\ref{m6}), and $R_\ep$, $T_\ep$, and $\Phi_\ep$
solve the \emph{acoustic equation}: \bFormula{r3} \ep \partial_t
(\alpha R_\ep + \beta T_\ep) +  \omega \Delta \Phi_\ep = 0, \eF
\bFormula{r4} \ep \partial_t \Grad \Phi_\ep + \Grad (\alpha R_\ep
+ \beta T_\ep) = 0, \eF with the initial data determined by
\bFormula{r4b} R_\ep (0, \cdot) = R_{0, \ep},\ T_\ep (0, \cdot) =
T_{0, \ep},\ \Phi_{\ep}(0, \cdot) = \Phi_{0, \ep}, \eF \noindent
where we have set
\[
\alpha = \frac{1}{\Ov{\vr}} \frac{\partial p (\Ov{\vr}, \Ov{\vt})}{\partial \vr}, \
\beta = \frac{1}{\Ov{\vr}} \frac{\partial p (\Ov{\vr}, \Ov{\vt})}{\partial \vt},\
\omega = \Ov{\vr} \left( \alpha + \frac{\beta^2}{\delta} \right).
\]

Noting that the functions $R_\ep$, $T_\ep$ are not uniquely determined by (\ref{r3} - \ref{r4b}),
we introduce the \emph{transport equation}
\bFormula{r4a}
\partial_t (\delta T_\ep - \beta R_\ep) + \vc{U}_\ep \cdot \Grad (\delta T_\ep - \beta R_\ep) +  (\delta T_\ep - \beta R_\ep) \Div \vc{U}_\ep = 0,
\eF
with
\[
\delta = \Ov{\vr} \frac{\partial s (\Ov{\vr}, \Ov{\vt})}{\partial \vt},
\]
where the initial data are determined by (\ref{r4b}). Equation (\ref{r4a}) is nothing other than a convenient linearization of the entropy
balance (\ref{i3}). Now, the system of equations (\ref{r3}), (\ref{r4}), (\ref{r4a}) is well-posed.

\subsection{Data regularization}

Our goal is to apply a Gronwall-type argument to the relative entropy inequality (\ref{r1}) to deduce the strong convergence to the limit system claimed in Theorem \ref{Tm1}. To this end, we choose the initial data
\bFormula{d1}
R_{0,\ep} = R_{0, \ep, \eta} = \chi_\eta * [ \psi_\eta \vr^{(1)}_{0, \ep}] ,\ T_{0,\ep} = T_{0, \ep, \eta} = \chi_\eta * [ \psi_\eta \vt^{(1)}_{0, \ep}], \ \eta > 0
\eF
where $\{ \chi_\eta (x) \}_{\eta > 0}$ is a family of regularizing kernels, and $\psi_\eta \in \DC(R^3)$ are the standard cut-off functions
$\psi_\eta \nearrow 1$. Similarly,
\bFormula{d2}
\Phi_{0,\ep} = \Phi_{0, \ep, \eta} = \chi_\eta * \Big[ \psi_\eta \Delta^{-1} \Div [\vu_{0, \ep}] \Big], \ \mbox{with}\
\Grad \Delta^{-1} \Div [\vu_{0,\ep}] \equiv \vc{H}^\perp [ \vu_{0,\ep} ].
\eF

To avoid excessive notation, we omit writing the parameter $\eta$ in the course of the limit passage $\ep \to 0$.

\section{Auxiliary  estimates}

We summarize the well known estimates for solutions of the auxiliary problems (\ref{r3}), (\ref{r4}), and (\ref{r4a}).

\label{a}

\subsection{Dispersive estimates}

The acoustic equation (\ref{r3} - \ref{r4b}) possesses a (unique)
smooth solution $\Phi_\ep$, $Z_\ep = \alpha R_\ep + \beta T_\ep$
satisfying the energy equality
\bFormula{a1} \left[ \left\| \Grad
\Phi_\ep (t, \cdot) \right\|^2_{W^{k,2}(R^3;R^3)} +
\frac{\delta}{\beta^2 + \alpha \delta} \left\| \alpha R_\ep(t,
\cdot) + \beta T_\ep (t, \cdot) \right\|^2_{W^{k,2}(R^3)}
\right]_{t = 0}^{t = \tau} = 0 \ \mbox{for all}\ { \tau \geq 0
},\ k = 0,1,2,\dots. \eF

In addition, we have the dispersive estimates
\bFormula{a2}
\left\| \Grad \Phi_\ep (t, \cdot) \right\|_{W^{k,q}(R^3;R^3)} + \left\| \alpha R_\ep(t, \cdot) + \beta T_\ep (t, \cdot) \right\|_{W^{k,q}(R^3)}
\eF
\[
\leq
c \left( 1 + \frac{t}{\ep} \right)^{-\left( \frac{1}{p} - \frac{1}{q} \right)} \left( \left\| \Grad \Phi_{0,\ep}  \right\|_{W^{d + k,p} (R^3;R^3)} + \left\| \alpha R_{0,\ep}  + \beta T_{0,\ep}  \right\|_{W^{d + k,p}(R^3)} \right),
\]
{for all $t \geq 0$,} where
\[
2 \leq q \leq \infty,\ \frac{1}{p} + \frac{1}{q} = 1,\ d > 3 \left( \frac{1}{p} - \frac{1}{q} \right) , \ k = 0,1,\dots,
\]
see Strichartz \cite{Strich}. Moreover, by virtue of the finite speed of propagation of acoustic waves,
the quantities $\Grad\Phi_\ep(t,\cdot)$ and $(\alpha R_\ep+\beta T_\ep)(t, \cdot)$ are compactly supported in $R^3$, see
(\ref{d1}), (\ref{d2}).

\subsection{Estimates for the transport equation}

The transport equation (\ref{r4a}) reads
\[
\partial_t (\delta T_\ep - \beta R_\ep) + (\vc{v} + \Grad \Phi_\ep)  \cdot \Grad (\delta T_\ep - \beta R_\ep) +  (\delta T_\ep - \beta R_\ep) \Delta \Phi_\ep = 0.
\]
In particular, we have \bFormula{a3} \Big[\intO{|\delta
T_\ep-\beta R_\ep|^2}\Big]_0^\tau= -
\int_0^\tau\intO{\Delta\Phi_\ep |\delta T_\ep-\beta R_\ep|^2}\
{\rm d}t, \eF and \bFormula{a4} \sup_{t \in [0,T]} \left\| \delta
T_\ep - \beta R_\ep \right\|_{W^{1,q}(\Omega)} \leq c(\eta, T)
\left\| \delta T_{0,\ep} - \beta R_{0,\ep}
\right\|_{W^{1,q}(\Omega)},\ 1 \leq q \leq \infty. \eF
{ Moreover, since the velocity of transport in the transport
equation is bounded and since $\Delta\Phi_\ep(t,\cdot)$ and the
initial data $\delta T_{0,\ep}-\beta R_{0,\ep}$ are compactly
supported, the solution $(\delta T_\ep-\beta R_\ep)(t,\cdot)$ is
as well compactly supported in $R^3$.}

\section{Convergence}
\label{c}

Fixing $\eta > 0$ our goal is perform the limit for $\ep \to 0$. This will be carried over in several steps.

\subsection{{ Viscous and heat conducting terms}}

{ We show that the dissipative terms related to viscosity and to
heat conductivity on the right-hand side of (\ref{r1}) are
negligible.} To this end, we write
\[
\ep^a \tn{S}(\vte, \Grad \vue) : \Grad \vc{U}_\ep = \ep^a \mu (\vte) \left( \Grad \vue + \Grad^t \vue - \frac{2}{3}
\Div \vue \tn{I} \right) : \Grad \vc{U}_\ep
\]
\[
= \ep^a \left[ \mu (\vte) \left( \Grad \vue + \Grad^t \vue - \frac{2}{3}
\Div \vue \tn{I} \right) \right]_{\rm ess} : \Grad \vc{U}_\ep +  \ep^a \left[ \mu (\vte) \left( \Grad \vue + \Grad^t \vue - \frac{2}{3}
\Div \vue \tn{I} \right) \right]_{\rm res} : \Grad \vc{U}_\ep,
\]
where
\[
\ep^a \int_{R^3} \left| \left[ \mu (\vte) \left( \Grad \vue + \Grad^t \vue - \frac{2}{3}
\Div \vue \tn{I} \right) \right]_{\rm ess} : \Grad \vc{U}_\ep \right| \ \dx \leq
\]
\[
\ep^{a/2} \left\|  \ep^{a/2} \left( \Grad \vue + \Grad^t \vue - \frac{2}{3}
\Div \vue \tn{I} \right)  \right\|_{L^2(R^3; R^{3 \times 3})} \left\| \Grad \vc{U}_{\ep} \right\|_{L^2(R^3;R^3)}  ;
\]
whence, by virtue of (\ref{b6}), (\ref{a1}),
\[
\ep^a \left[ \mu (\vte) \left( \Grad \vue + \Grad^t \vue - \frac{2}{3}
\Div \vue \tn{I} \right) \right]_{\rm ess} : \Grad \vc{U}_\ep \to 0 \ \mbox{in} \ L^2((0,T) \times \Omega)
\ \mbox{as} \ \ep \to 0.
\]

Similarly, in accordance with (\ref{b5}), (\ref{b6}), and hypothesis (\ref{i8}),
\[
\ep^a \left[ \mu (\vte) \left( \Grad \vue + \Grad^t \vue - \frac{2}{3}
\Div \vue \tn{I} \right) \right]_{\rm res} : \Grad \vc{U}_\ep
\]
\[
= \ep^a \sqrt{ [\vte]_{\rm res} } \sqrt{ [ \mu (\vte) ]_{\rm res}
} \sqrt{ \frac{\mu (\vte)}{\vte}}  \left( \Grad \vue + \Grad^t
\vue - \frac{2}{3} \Div \vue \tn{I} \right)  : \Grad \vc{U}_\ep
\to 0 \ \mbox{in}\ L^2(0,T; L^{4/3}(\Omega;R^3)) \ \mbox{as}\ \ep
\to 0.
\]

Next, we have
\[
\ep^{b - 2} \frac{ \vc{q}(\vte, \Grad \vte) \cdot \Grad \Theta_\ep }{\vte} = - \ep^{b-1} \frac{ \kappa(\vte) \Grad (\vte - \Ov{\vt})}{\vte} \cdot
\Grad T_\ep
\]
\[
= - \ep^{b/2}  \left[ \ep^{(b-2)/2} \frac{ \kappa(\vte) }{\vte} \Grad (\vte - \Ov{\vt}) \right]_{\rm ess} \cdot \Grad T_\ep
- \ep^{b/2} \left[ \ep^{(b-2)/2} \frac{ \kappa(\vte) }{\vte} \Grad (\vte - \Ov{\vt}) \right]_{\rm res} \cdot \Grad T_\ep,
\]
where, as a consequence of (\ref{b7}), (\ref{a2}), (\ref{a4}),
\[
\ep^{b/2}  \left[ \ep^{(b-2)/2} \frac{ \kappa(\vte) }{\vte} \Grad (\vte - \Ov{\vt}) \right]_{\rm ess} \cdot \Grad T_\ep \to 0
\ \mbox{in}\ L^2((0,T) \times R^3) \ \mbox{as}\ \ep \to 0.
\]

Moreover, in accordance with hypothesis (\ref{i9}),
\[
\ep^{b/2} \left| \left[ \ep^{(b-2)/2} \frac{ \kappa(\vte) }{\vte} \Grad (\vte - \Ov{\vt}) \right]_{\rm res} \cdot \Grad T_\ep \right|
\]
\[
\leq c \ep^{b/2} \left( \left| \ep^{(b-2)/2} \Grad \Big( \log(\vte) - \log(\Ov{\vt} ) \Big) \right| +
\left| \ep^{(b-2)/2}  \left[ \vte \right]_{\rm res}^2 \Grad \Big( \vte - \Ov{\vt}  \Big) \right|
 \right) \left| \Grad T_\ep \right|,
\]
where, by virtue of (\ref{b5}), (\ref{b7}), (\ref{a2}), and (\ref{a4}),  the right-hand side tends to zero in $L^1((0,T) \times R^3)$.

Thus (\ref{r1}) reduces to

\greybox{

\bFormula{r9}
\left[ \mathcal{E}_\ep \left( \vre, \vte, \vue \Big| r_\ep , \Theta_\ep , \vc{U}_\ep \right) \right]_{t = 0}^\tau
\eF
\[
\leq \int_0^\tau \intO{  \vre \Big( \partial_t \vc{U}_\ep + \vue
\cdot\Grad \vc{U}_\ep \Big) (\vc{U}_\ep - \vue)    } \ \dt
\]
$$
+\frac
1{\ep^2}\int_0^\tau\intO{\Big[\Big(p(r_\ep,\Theta_\ep)-p(\vre,\vte)\Big){\rm
div}\vc U_\ep +\frac\vre {r_\ep}(\vc U_\ep-\vue)\cdot\Grad
p(r_\ep,\Theta_\ep)\Big]}{\rm d}t
$$
\[
- \frac{1}{\ep} \int_0^\tau \intO{ \left( \vre \Big( s(\vre,\vte) - s(r_\ep, \Theta_\ep) \Big) \partial_t T_\ep + \vre \Big( s(\vre,\vte) - s(r_\ep, \Theta_\ep) \Big) \vue \cdot \Grad T_\ep \right) } \ \dt
\]
\[
+ \frac{1}{\ep^2} \int_0^\tau \intO{ \frac{r_\ep-\vre}{r_\ep}\Big(
\partial_t p(r_\ep , \Theta_\ep ) +\vc U_\ep \cdot
\Grad p(r_\ep , \Theta_\ep ) \Big) } \ \dt + \chi_\ep (\tau, \eta)
\]
with
\[
\chi_\ep (\cdot, \eta) \to 0 \ \mbox{in}\ C[0,T] \ \mbox{as} \ \ep \to 0 \ \mbox{for any fixed}\ \eta > 0.
\]

}

\subsection{Velocity dependent terms}

Our goal is to handle the integral
\[
\int_0^\tau \intO{ \Big[  \vre (\vc{U}_\ep - \vue ) \cdot
\partial_t \vc{U}_\ep + \vre (\vc{U}_\ep - \vue) \otimes \vue : \Grad \vc{U}_\ep \Big] } \ \dt =
\]
\[
\int_0^\tau \intO{ \Big[  \vre (\vc{U}_\ep - \vue ) \cdot
\partial_t \vc{U}_\ep + \vre (\vc{U}_\ep - \vue) \otimes \vc{U}_\ep : \Grad \vc{U}_\ep \Big] } \ \dt + \int_0^\tau \intO{ \vre (\vc{U}_\ep - \vue ) \otimes
(\vue - \vc{U}_\ep ) : \Grad \vc{U}_\ep } \ \dt,
\]
{where the second term on the right-hand side is bounded by
$$
\int_0^\tau \Big\|\Grad\vc
U_\ep\Big\|_{L^\infty(\Omega;R^{3\times
3})}\mathcal{E}_\ep \left( \vre, \vte, \vue \Big| r_\ep ,
\Theta_\ep , \vc{U}_\ep \right){\rm d}t
\le
\int_0^\tau \big(c+\chi_\ep(t,\eta)\Big)\mathcal{E}_\ep \left(
\vre, \vte, \vue \Big| r_\ep , \Theta_\ep , \vc{U}_\ep \right){\rm
d}t,
$$ with $c$ independent of $\ep$, $\eta$,
and, by virtue of
(\ref{euler}) and (\ref{a2}), $\chi_\ep (\cdot,\eta)\to 0$ in  $C[0,T]$ as $\ep \to 0$.}


On the other hand,
\[
\int_0^\tau \intO{ \Big[  \vre (\vc{U}_\ep - \vue ) \cdot
\partial_t \vc{U}_\ep + \vre (\vc{U}_\ep - \vue) \otimes \vc{U}_\ep : \Grad \vc{U}_\ep \Big] } \ \dt
\]
\[
= \int_0^\tau \intO{ \vre (\vc{U}_\ep - \vue) \cdot \Big( \partial_t \vc{v} + \vc{v}
\cdot \Grad \vc{v} \Big) } \ \dt + \int_0^\tau \intO{ \vre (\vc{U}_\ep - \vue)
\cdot \partial_t \Grad \Phi_\ep } \ \dt
\]
\[
+ \int_0^\tau \intO{ \vre (\vc{U}_\ep - \vue ) \otimes
\Grad \Phi_\ep : \Grad \vc{v}} + \int_0^\tau \intO{ \vre (\vc{U}_\ep - \vue ) \otimes
\vc{v} : \Grad^2 \Phi_\ep } \ \dt
\]
\[
+\frac{1}{2} \int_0^\tau \intO{ \vre (\vc{U}_\ep - \vue) \cdot \Grad |\Grad \Phi_\ep |^2 } \ \dt.
\]

In view of the uniform bounds (\ref{euler}), (\ref{b3} - \ref{b5}), and the {dispersive estimates} stated in (\ref{a2}),
the last three integrals tend to zero for $\ep \to 0$, uniformly with respect to $\tau$. Accordingly, we focus on the first two terms, where the former reads
\[
 \int_0^\tau \intO{ \vre (\vc{U}_\ep - \vue) \cdot \Big( \partial_t \vc{v} + \vc{v}
\cdot \Grad \vc{v} \Big) } \ \dt = - \int_0^\tau \intO{ \vre (\vc{U}_\ep - \vue) \cdot \Grad \Pi } \  \dt
\]
\[
=
\int_0^\tau \intO{ \vre \vue \cdot \Grad \Pi } \ \dt - \int_0^\tau \intO{
\vre (\vc{v} + \Grad \Phi_\ep ) \cdot \Grad \Pi } \ \dt.
\]

As a consequence of the estimates (\ref{b3} - \ref{b5}), we get
\bFormula{conv1}
\vre \vue \to \Ov{ \vr \vu } \ \mbox{weakly-(*) in}\ L^\infty(0,T; L^2 + L^{5/4}(R^3;R^3)),
\eF
where, thanks to the continuity equation (\ref{v1}),
\bFormula{conv2}
\Div (\Ov{ \vr \vu }) = 0,
\eF
in particular,
\[
\int_0^\tau \intO{ \vre \vue \Grad \Pi } \ \dt \to 0 \ \mbox{in}\ L^q(0,T) \ \mbox{for any}\ 1 \leq q < \infty.
\]

Next, we have
\[
\int_0^\tau \intO{ \vre \left( \vc{v} + \Grad \Phi_\ep \right) \cdot \Grad \Pi }
\ \dt
\]
\[
= \ep \int_0^\tau \intO{ \frac{\vre - \Ov{\vr}}{\ep} \vc{v} \cdot \Grad \Pi } \ \dt +
\ep \int_0^\tau \intO{ \frac{\vre - \Ov{\vr}}{\ep}   \Grad \Phi_\ep \cdot \Grad \Pi  } \ \dt +
\int_0^\tau \intO{ \Ov{\vr}   \Grad \Phi_\ep \cdot \Grad \Pi  } \ \dt,
\]
where the first two integrals vanish in the limit $\ep \to 0$, while
\[
\int_0^\tau \intO{ \Grad \Phi_\ep \cdot \Grad \Pi  } \ \dt = - \int_0^\tau \intO{ \Delta \Phi_\ep \cdot \Pi  } \ \dt =
\frac{\ep \Ov{\vr}}{\omega} \int_0^\tau \intO{ \partial_t (\alpha R_\ep + \beta T_\ep ) \Pi } \ \dt
\]
\[
{ =}\frac{\ep \Ov{\vr}}{\omega} \left[ \intO{  (\alpha R_\ep +
\beta T_\ep ) \Pi } \right]_{t = 0}^{\tau} - \frac{\ep
\Ov{\vr}}{\omega} \int_0^\tau \intO{ (\alpha R_\ep + \beta T_\ep )
\partial_t \Pi } \ \dt = \chi_\ep (\tau, \eta),
\]
where here and hereafter, the symbol $\chi_\ep(\tau,\eta)$ denotes a generic function satisfying
\[
\chi_\ep(\cdot, \eta) \to 0 \ \mbox{in}\ L^1(0,T) \ \mbox{as}\ \ep \to 0 \ \mbox{for any fixed}\ \eta > 0.
\]

Thus, it remains to handle
\[
\int_0^\tau \intO{ \vre (\vc{U}_\ep - \vue ) \cdot \partial_t
\Grad \Phi_\ep } \ \dt
\]
\[
 = - \int_0^\tau \intO{ \vre \vue \cdot \partial_t \Grad \Phi_\ep } \ \dt +
\int_0^\tau \intO{ \vre \vc{v} \cdot \partial_t \Grad \Phi_\ep } \ \dt +
\frac{1}{2} \int_0^\tau \intO{ \vre \partial_t | \Grad \Phi_\ep |^2 } \ \dt,
\]
where, in accordance with (\ref{b4}), (\ref{b5}), and the dispersive estimates (\ref{a2}),
\[
\int_0^\tau \intO{ \vre \vc{v} \cdot \partial_t \Grad \Phi_\ep } \ \dt =
\int_0^\tau \intO{ (\vre - \Ov{\vr})  \vc{v} \cdot \partial_t \Grad \Phi_\ep } \ \dt
\]
\[
 = - \int_0^\tau \intO{ \frac{\vre - \Ov{\vr}}{\ep}  \vc{v} \cdot \Grad \left( \alpha R_\ep + \beta T_\ep \right) } \ \dt = \chi_\ep (\eta, \tau)
\]
while
\[
\frac{1}{2} \int_0^\tau \intO{ \vre \partial_t | \Grad \Phi_\ep |^2 } \ \dt = \frac{\ep}{2} \int_0^\tau \intO{ \frac{\vre - \Ov{\vr}}{\ep} \partial_t | \Grad \Phi_\ep |^2 } \ \dt +
\frac{1}{2} \int_0^\tau \intO{ \Ov{\vr} \partial_t | \Grad \Phi_\ep |^2 } \ \dt .
\]
Finally, using (\ref{r4}), we get
\[
\frac{\ep}{2} \int_0^\tau \intO{ \frac{\vre - \Ov{\vr}}{\ep} \partial_t | \Grad \Phi_\ep |^2 } \ \dt = -
\int_0^\tau \intO{ \frac{\vre - \Ov{\vr}}{\ep} \Grad \Phi_\ep \cdot \Grad (\alpha R_\ep + \beta T_\ep) } \ \dt,
\]
where, by virtue of the dispersive estimates (\ref{a2}), the last integral tends to zero.

Thus relation (\ref{r9}) reduces to

\greybox{

\bFormula{r10}
\left[ \mathcal{E}_\ep \left( \vre, \vte, \vue \Big| r_\ep , \Theta_\ep , \vc{U}_\ep \right) \right]_{t = 0}^\tau
\eF
\[
\leq \left[ \intO{ \Ov{\vr} \frac{1}{2}|\Grad \Phi_\ep |^2 } \right]_{t = 0}^{t = \tau} -
\int_0^\tau \intO{ \vre \vue \cdot \partial_t \Grad \Phi_\ep } \ \dt
\]
\[
- \frac{1}{\ep} \int_0^\tau \intO{ \left[ \vre \Big( s(\vre,\vte) - s(r_\ep, \Theta_\ep) \Big) \partial_t T_\ep + \vre \Big( s(\vre,\vte) - s(r_\ep, \Theta_\ep) \Big) \vue \cdot \Grad T_\ep  \right] } \ \dt
\]
\[
+ \frac{1}{\ep^2} \int_0^\tau \intO{ \left( ( r_\ep - \vre ) \frac{1}{r_\ep } \partial_t p(r_\ep , \Theta_\ep ) - \frac{\vre}{r_\ep} \vue \cdot
\Grad p(r_\ep , \Theta_\ep ) \right) } \ \dt - \frac{1}{\ep^2} \int_0^\tau \intO{ \Big( p(\vre, \vte) - p(\Ov{\vr}, \Ov{\vt}) \Big)
\Delta \Phi_\ep } \ \dt
\]
\[
+ \int_0^\tau \left(c + \chi^1_{\ep} (t, \eta) \right) \mathcal{E}_\ep \left( \vre, \vte, \vue \Big| r_\ep , \Theta_\ep , \vc{U}_\ep \right) \ \dt   + \chi^2_\ep
(\tau, \eta),
\]
where
\[
\chi^i_{\ep}(\cdot, \eta) \to 0 \ \mbox{in} \ L^1(0,T) \ \mbox{as}\ \ep \to 0 \ \mbox{for any fixed}\ \eta > 0, \ i=1,2,
\]

}

\noindent and where we have used the identity
\[
\intO{ \left[ \Big( p(r_\ep, \Theta_\ep) - p(\vre, \vte) \Big) \Div \vc{U}_\ep + \left( 1 - \frac{\vre}{r_\ep} \right) \vc{U}_\ep
\cdot \Grad p(r_\ep, \Theta_\ep) + \frac{\vre}{r_\ep} (\vc{U}_\ep - \vue) \cdot \Grad p(r_\ep, \Theta_\ep)  \right] }
\]
\[
= - \intO{ \Big( p(\vre, \vte) - p(\Ov{\vr}, \Ov{\vt}) \Big) \Delta \Phi_\ep}  - \intO{ \frac{\vre}{r_\ep} \vue \cdot \Grad p(r_\ep, \Theta_\ep) }.
\]
Recall that $\Grad \Phi_\ep(t, \cdot)$ is compactly supported and $\Div \vc{v} = 0$, which justifies the by-parts integration used in the above.

\subsection{Pressure terms}

We write
\[
\frac{1}{\ep^2} \vre \vue\cdot \frac{1}{r_\ep} \Grad p(r_\ep,
\Theta_\ep) = \frac{1}{\ep} \vre \vue\cdot \frac{1}{r_\ep} \left(
\frac{ \partial p(r_\ep, T_\ep) }{\partial \vr} \Grad R_\ep +
\frac{ \partial p(r_\ep, T_\ep) }{\partial \vt} \Grad T_\ep
\right)
\]
\[
= \frac{1}{\ep} \vre \vue \cdot\frac{1}{r_\ep} \left[ \left(
\frac{
\partial p(r_\ep, T_\ep) }{\partial \vr} - \frac{ \partial
p(\Ov{\vr}, \Ov{\vt} ) }{\partial \vr} \right)
 \Grad R_\ep + \left( \frac{ \partial p(r_\ep, T_\ep) }{\partial \vt} - \frac{ \partial p(\Ov{\vr}, \Ov{\vt} ) }{\partial \vt} \right) \Grad T_\ep \right]
+ \frac{1}{\ep} \vre \vue\cdot \frac{\Ov{\vr}}{r_\ep} \Grad \left(
\alpha R_\ep + \beta T_\ep \right)
\]
\[
= \frac{1}{\ep} \vre \vue\cdot \frac{1}{r_\ep} \left[ \left(
\frac{
\partial p(r_\ep, T_\ep) }{\partial \vr} - \frac{ \partial
p(\Ov{\vr}, \Ov{\vt} ) }{\partial \vr} \right)
 \Grad R_\ep + \left( \frac{ \partial p(r_\ep, T_\ep) }{\partial \vt} - \frac{ \partial p(\Ov{\vr}, \Ov{\vt} ) }{\partial \vt} \right) \Grad T_\ep \right]
\]
\[
+  \frac{1}{\ep} \vre \vue \cdot \Grad \left( \alpha R_\ep + \beta
T_\ep \right) +\frac{1}{\ep} \vre \vue\cdot \left(
\frac{\Ov{\vr}}{r_\ep} - 1 \right)  \Grad \left( \alpha R_\ep +
\beta T_\ep \right),
\]
where, by virtue of (\ref{b3}), (\ref{b4}), and the dispersive estimates (\ref{a2}),
\[
\frac{1}{\ep} \vre \vue\cdot \left( \frac{\Ov{\vr}}{r_\ep} - 1
\right) \Grad \left( \alpha R_\ep + \beta T_\ep \right) \to 0 \
\mbox{in}\ L^q(0,T; L^2 + L^{5/4} (R^3;R^3)) ,\ 1 \leq q < \infty
\ \mbox{for}\ \ep \to 0,
\]
while, in accordance with (\ref{r4}),
\[
\frac{1}{\ep} \vre \vue \cdot \Grad \left( \alpha R_\ep + \beta T_\ep \right) = -\vre \vue \cdot \partial_t \Grad \Phi_\ep.
\]

Finally, using the Taylor expansion formula, we obtain
\[
\int_0^\tau \intO{ \frac{1}{\ep}
\vre \vue \cdot \frac{1}{r_\ep} \left[ \left( \frac{ \partial p(r_\ep, T_\ep) }{\partial \vr} - \frac{ \partial p(\Ov{\vr}, \Ov{\vt} ) }{\partial \vr} \right)
 \Grad R_\ep + \left( \frac{ \partial p(r_\ep, T_\ep) }{\partial \vt} - \frac{ \partial p(\Ov{\vr}, \Ov{\vt} ) }{\partial \vt} \right) \Grad T_\ep \right] } \ \dt
\]
\[
= \int_0^\tau \intO{ \vre \vue \cdot \left[ \frac{1}{2} \frac{ \partial^2 p(\Ov{\vr}, \Ov{\vt} ) }{\partial \vr^2} \Grad R^2_\ep+ \frac{\partial^2 p(\Ov{\vr}, \Ov{\vt})}{\partial \vr \partial \vt } \Grad (R_\ep T_\ep)
 + \frac{1}{2} \frac{ \partial^2 p(\Ov{\vr}, \Ov{\vt} ) }{\partial \vt^2}  \Grad T^2_\ep \right]} \ \dt + \chi_\ep(\tau, \eta);
\]
where, furthermore, as $\vre \vue$ satisfies (\ref{conv1}), (\ref{conv2}),
\[
\int_0^\tau \intO{
\vre \vue \cdot \left[ \frac{1}{2} \frac{ \partial^2 p(\Ov{\vr}, \Ov{\vt} ) }{\partial \vr^2} \Grad R^2_\ep+ \frac{\partial^2 p(\Ov{\vr}, \Ov{\vt})}{\partial \vr \partial \vt } \Grad (R_\ep T_\ep)
 + \frac{1}{2} \frac{ \partial^2 p(\Ov{\vr}, \Ov{\vt} ) }{\partial \vt^2}  \Grad T^2_\ep \right] } \ \dt \to 0
 \ \mbox{in}\ L^1(0,T).
\]

Consequently, we may infer that (\ref{r10}) reduces to

\greybox{

\bFormula{r11}
\left[ \mathcal{E}_\ep \left( \vre, \vte, \vue \Big| r_\ep , \Theta_\ep , \vc{U}_\ep \right) \right]_{t = 0}^\tau
\leq \left[ \intO{ \Ov{\vr} \frac{1}{2}|\Grad \Phi_\ep |^2 } \right]_{t = 0}^{t = \tau}
\eF
\[
- \frac{1}{\ep} \int_0^\tau \intO{ \left[ \vre \Big( s(\vre,\vte) - s(r_\ep, \Theta_\ep) \Big) \partial_t T_\ep + \vre \Big( s(\vre,\vte) - s(r_\ep, \Theta_\ep) \Big) \vue \cdot \Grad T_\ep \right] } \ \dt
\]
\[
+ \frac{1}{\ep^2} \int_0^\tau \intO{\frac{r_\ep - \vre}
{r_\ep} \partial_t p(r_\ep , \Theta_\ep )  } \ \dt
- \frac{1}{\ep^2} \int_0^\tau \intO{ \Big(p(\vre,
\vte)-p(\overline\vr,\overline\vt)\Big) \Delta \Phi_\ep } \ \dt
\]
\[
+ \int_0^\tau \left(c + \chi^1_{\ep} (t, \eta) \right) \mathcal{E}_\ep \left( \vre, \vte, \vue \Big| r_\ep , \Theta_\ep , \vc{U}_\ep \right) \ \dt   + \chi^2_\ep
(\tau, \eta),
\]
where
\[
\chi^i_{\ep}(\cdot, \eta) \to 0 \ \mbox{in} \ L^1(0,T) \ \mbox{as}\ \ep \to 0 \ \mbox{for any fixed}\ \eta > 0, \ i=1,2.
\]

}

\subsection{Replacing velocity in the convective term}

Our next goal is to ``replace'' $\vue$ by $\vc{U}_\ep$ in the remaining convective term in (\ref{r11}). To this end,
we write
\[
\int_0^\tau \intO{ \vre \frac{ s(\vre, \vte) - s(r_\ep, \Theta_\ep ) }{ \ep } \vue \cdot \Grad T_\ep } \ \dt
\]
\[
=
\int_0^\tau \intO{ \vre \frac{ s(\vre, \vte) - s(r_\ep , \Theta_\ep ) }{ \ep } \vc{U}_\ep  \cdot \Grad T_\ep } \ \dt +
\int_0^\tau \intO{ \vre \frac{ s(\vre, \vte) - s(r_\ep, \Theta_\ep ) }{ \ep } (\vue - \vc{U}_\ep)  \cdot \Grad T_\ep } \ \dt,
\]
where
\[
\int_0^\tau \intO{ \vre \frac{ s(\vre, \vte) - s(r_\ep, \Theta_\ep ) }{ \ep } (\vue - \vc{U}_\ep)  \cdot \Grad T_\ep } \ \dt\]
\[
= \int_0^\tau \intO{ \vre \left[ \frac{ s(\vre, \vte) - s(r_\ep, \Theta_\ep ) }{ \ep } \right]_{\rm ess} (\vue - \vc{U}_\ep)  \cdot \Grad T_\ep } \ \dt
\]
\[
+ \int_0^\tau \intO{ \vre \left[ \frac{ s(\vre, \vte) - s(r_\ep , \Theta_\ep ) }{ \ep } \right]_{\rm res} (\vue - \vc{U}_\ep)  \cdot \Grad T_\ep } \ \dt.
\]

Next, we get
\[
\left| \int_0^\tau \intO{ \vre \left[ \frac{ s(\vre, \vte) - s(r_\ep , \Theta_\ep ) }{ \ep } \right]_{\rm ess} (\vue - \vc{U}_\ep)  \cdot \Grad T_\ep } \ \dt \right|
\]
\[
\leq c \int_0^\tau \| \Grad T_\ep(t, \cdot) \|_{L^\infty(R^3;R^3)} \intO{ \left( \vre | \vue - \vc{U}_\ep |^2
+ \left| \left[ \frac{\vre - r_\ep}{\ep} \right]_{\rm ess} \right|^2 + \left| \left[ \frac{\vte - \Theta_\ep}{\ep} \right]_{\rm ess} \right|^2 \right)
}\ \dt;
\]
whence this term can be ``absorbed'' by means of Gronwall argument.

As for the residual component, we have to control the most difficult term
$[\vre s(\vre, \vte)]_{\rm res} \vue$. To begin, the hypotheses (\ref{i12} - \ref{i16}) imply that
\[
\vr | s(\vr, \vt) | \leq c\left( \vt^3 + \vr |\log(\vr)| + \vr [\log(\vt)]^+ \right).
\]
Consequently, by virtue of the estimates (\ref{b5}), (\ref{b6}),
\[
\left\| \left[ \vte^3 \right]_{\rm res} \vue \right\|_{L^1(R^3;R^3)} \leq \ep^{-a/2} \left\| \left[ \vte^3 \right]_{\rm res} \right\|_{L^{6/5}(R^3)}
\| \ep^{a/2} \vue \|_{W^{1,2}(R^3;R^3)}
\]
\[
c_2 \leq \ep^{\left(\frac{5}{3} - \frac{a}{2} \right)} \| \ep^{a/2} \vue \|_{W^{1,2}(R^3;R^3)} \to 0 \ \mbox{in}\ L^2(0,T) \ \mbox{whenever}\ 0 < a < \frac{10}{3}.
\]

Estimating the remaining integrals in a similar way, we
can rewrite inequality (\ref{r11}) in the form

\greybox{

\bFormula{r11+}
\left[ \mathcal{E}_\ep \left( \vre, \vte, \vue
\Big| r_\ep , \Theta_\ep , \vc{U}_\ep \right) \right]_{t = 0}^\tau
\leq \left[ \intO{ \Ov{\vr} \frac{1}{2}|\Grad \Phi_\ep |^2 }
\right]_{t = 0}^{t = \tau} \eF
\[
- \frac{1}{\ep} \int_0^\tau \intO{ \left[ \vre \Big( s(\vre,\vte)
- s(r_\ep, \Theta_\ep) \Big) \partial_t T_\ep + \vre \Big(
s(\vre,\vte) - s(r_\ep, \Theta_\ep) \Big) \vc U_\ep \cdot \Grad
T_\ep \right] } \ \dt
\]
\[
+ \frac{1}{\ep^2} \int_0^\tau \intO{ \frac{r_\ep - \vre}
{r_\ep}\partial_t p(r_\ep , \Theta_\ep )  } \ \dt -
\frac{1}{\ep^2} \int_0^\tau \intO{ \Big(p(\vre,
\vte)-p(\overline\vr,\overline\vt)\Big) \Delta \Phi_\ep } \ \dt
\]
\[
+ \int_0^\tau c\left(1 + \chi^1_{\ep} (t, \eta) + \| \Grad T_\ep (t, \cdot) \|_{L^\infty(R^3;R^3)} \right) \mathcal{E}_\ep \left( \vre, \vte, \vue \Big| r_\ep , \Theta_\ep , \vc{U}_\ep \right) \ \dt   + \chi^2_\ep
(\tau, \eta),
\]
where
\[
\chi^i_{\ep}(\cdot, \eta) \to 0 \ \mbox{in} \ L^1(0,T) \ \mbox{as}\ \ep \to 0 \ \mbox{for any fixed}\ \eta > 0, \ i=1,2.
\]

}

\subsection{Entropy and pressure}

In order to handle the remaining integrals in (\ref{r11+}), we first show that all terms can be replaced by their linearization
at $\Ov{\vr}$, $\Ov{\vt}$. To this end, we first observe that we may neglect the ``residual part'' of all integrals. Indeed,
\[
\frac{1}{\ep} \int_0^\tau \intO{ \left[ \vre \Big( s(\vre,\vte)
- s(r_\ep, \Theta_\ep) \Big) \right]_{\rm res} \partial_t T_\ep } \ \dt =
\frac{1}{\ep^2} \int_0^\tau \intO{ \left[ \vre \Big( s(\vre,\vte)
- s(r_\ep, \Theta_\ep) \Big) \right]_{\rm res} \ep {\partial_t T_\ep}} \ \dt,
\]
where, by virtue of the estimates (\ref{a2} - \ref{a4}), the equations (\ref{r3} - \ref{r4a}), and the identities,
\bFormula{ident}
(\beta^2+\alpha\delta) T=\beta(\alpha R+\beta T)+\alpha(\delta
T-\beta R),
\
(\beta^2+\alpha\delta) R=\delta(\alpha R+\beta T)-\beta(\delta
T-\beta R),
\eF
we get
\bFormula{pom1}
\sup_{t \in [0,T]} \ep \| \partial_t R_\ep (t, \cdot) \|_{L^\infty(R^3)} ,\
\sup_{t \in [0,T]} \ep \| \partial_t T_\ep (t, \cdot) \|_{L^\infty(R^3)} \leq c(\eta),
\eF
\bFormula{pom2}
\ep \| \partial_t R_\ep (t, \cdot) \|_{L^\infty(R^3)} \to 0 , \ \ep \| \partial_t T_\ep (t, \cdot) \|_{L^\infty(R^3)} \to 0\ \mbox{for any} \ t > 0,
\eF
while, in accordance with (\ref{b5}),
\[
{\rm ess} \sup_{t \in (0,T)} \intO{ \vre \left[ s(\vre,\vte)
- s(r_\ep, \Theta_\ep) \Big) \right]_{\rm res} } \leq \ep^2 c.
\]
A similar treatment can be applied to the integrals
\[
\frac{1}{\ep} \int_0^\tau \intO{
\vre \left[
s(\vre,\vte) - s(r_\ep, \Theta_\ep) \right]_{\rm res} \vc U_\ep \cdot \Grad
T_\ep } \ \dt
\ \mbox{and}\
\frac{1}{\ep^2} \int_0^\tau \intO{
\left[ p(\vre,
\vte)-p(\overline\vr,\overline\vt)\right]_{\rm res} \Delta \Phi_\ep } \ \dt.
\]

Finally,
\[
\frac{1}{\ep^2} \int_0^\tau \intO{ \left[ \frac{r_\ep - \vre}
{r_\ep} \right]_{\rm res} \partial_t p(r_\ep , \Theta_\ep )  } \ \dt
\]
\[
= \frac{1}{\ep^2} \int_0^\tau \intO{ \left[ \frac{r_\ep - \vre}
{r_\ep} \right]_{\rm res} \left( \frac{\partial p(r_\ep, \Theta_\ep)}{\partial \vr} \ep \partial_t  R_\ep + \frac{\partial p(r_\ep, \Theta_\ep)}{\partial \vt} \ep \partial_t  T_\ep
\right)  } \ \dt;
\]
whence (\ref{pom1}), (\ref{pom2}) yield the desired conclusion.

Since all remaining integrals in (\ref{r11+}) can be reduced to their ``essential component'' it is easy to check that
\begin{equation}\label{*1}
- \frac{1}{\ep} \int_0^\tau \intO{ \left[ \vre \Big( s(\vre,\vte)
- s(r_\ep, \Theta_\ep) \Big) \partial_t T_\ep + \vre \Big(
s(\vre,\vte) - s(r_\ep, \Theta_\ep) \Big) \vc U_\ep \cdot \Grad
T_\ep \right] } \ \dt
\end{equation}
$$
- \frac{1}{\ep^2} \int_0^\tau \intO{ \frac{\vr_\ep - r_\ep}
{r_\ep}\partial_t p(r_\ep , \Theta_\ep )  } \ \dt -\frac 1{\ep^2}
\int_0^\tau
\intO{\Big(p(\vre,\vte)-p(\overline\vr,\overline\vt)\Big)\Delta\Phi_\ep}\ {\rm
d}t
$$
$$
= -\int_0^\tau\intO{\Big(\delta\frac{\vte-\Theta_\ep}\ep
-\beta\frac{\vr_\ep-r_\ep}\ep\Big)\Big(\partial_t T_\ep + \vc
U_\ep\cdot\Grad T_\ep\Big)}\ {\rm d}t
$$
$$
-\int_0^\tau\intO{\frac{\vre- r_\ep}\ep\partial_t\Big(\alpha
R_\ep+\beta T_\ep\Big)}{\rm d}t +\int_0^\tau
\intO{\frac{\delta}{\beta^2+\alpha\delta}\Big(\alpha\frac{\vre
-\overline\vr}\ep+\beta\frac{\vte-\overline\vt}\ep\Big)
\partial_t\Big(\alpha R_\ep+\beta T_\ep\Big)}\ {\rm d}t
+ \chi_\ep(\tau, \eta)
$$
$$
= \int_0^\tau\intO{\Big(\delta T_\ep-\beta R_\ep\Big)\partial_t
T_\ep}{\rm d} t + \int_0^\tau\intO{R_\ep\partial_t\Big(\alpha
R_\ep+\beta T_\ep\Big)}\ {\rm d}t
$$
$$
-\left[ \int_0^\tau\intO{\Big(\delta\frac{\vte-\overline\vt}\ep
-\beta\frac{\vr_\ep-\overline\vr}\ep\Big)\partial_t T_\ep}\ {\rm d}t
+
\int_0^\tau\intO{\Big(\frac{\beta^2}{\beta^2+\alpha\delta}\frac{\vre-\overline\vr}\ep
-\frac{\beta\delta}{\beta^2+\alpha\delta}\frac{\vte-\overline\vt}\ep\Big)
\partial_t\Big(\alpha
R_\ep+\beta T_\ep\Big)}\ {\rm d}t \right]
$$
$$
- \int_0^\tau\intO{\Big(\delta\frac{\vte-\Theta_\ep}\ep
-\beta\frac{\vr_\ep-r_\ep}\ep\Big)\vc U_\ep\cdot\Grad T_\ep} \
{\rm d}t + \chi_\ep(\tau, \eta),
$$
{ where we have used (\ref{A}--\ref{r4}).}

In the next step, we use the identities (\ref{ident})
to compute,
\begin{equation}\label{*2}
\int_0^\tau\intO{\Big(\delta T_\ep-\beta R_\ep\Big)\partial_t
T_\ep}\ {\rm d} t + \int_0^\tau\intO{R_\ep\partial_t\Big(\alpha
R_\ep+\beta T_\ep\Big)}\ {\rm d}t
\end{equation}
$$
=
\int_0^\tau\intO{\left[\frac\beta{\beta^2+\alpha\delta}\Big(\delta
T_\ep-\beta R_\ep\Big)\partial_t\Big(\alpha R_\ep+\beta T_\ep\Big)
+ \frac\alpha{\beta^2+\alpha\delta}\Big(\delta T_\ep-\beta
R_\ep\Big)\partial_t\Big(\delta T_\ep-\beta R_\ep\Big) \right.
$$
$$
\left.
+\frac\delta{\beta^2+\alpha\delta}\Big(\alpha R_\ep+\beta
T_\ep\Big)\partial_t\Big(\alpha R_\ep+\beta T_\ep\Big) -
\frac\beta{\beta^2+\alpha\delta}\Big(\delta T_\ep-\beta
R_\ep\Big)\partial_t\Big(\alpha R_\ep+\beta T_\ep\Big)\right]}\ {\rm
d}t
$$
$$
=\frac12\frac\delta{\beta^2+\alpha\delta} \left[\intO{|\alpha
R_\ep+\beta T_\ep|^2}\right]_0^\tau+
\frac12\frac\alpha{\beta^2+\alpha\delta} \left[\intO{|\delta
T_\ep-\beta R_\ep|^2}\right]_0^\tau.
$$

Similarly, we get
\begin{equation}\label{*3}
- \int_0^\tau\intO{\Big(\delta\frac{\vte-\overline\vt}\ep
-\beta\frac{\vr_\ep-\overline\vr}\ep\Big)\partial_t T_\ep}\ {\rm d}t
-
\int_0^\tau\intO{\Big(\frac{\beta^2}{\beta^2+\alpha\delta}\frac{\vre-\overline\vr}\ep
-\frac{\beta\delta}{\beta^2+\alpha\delta}\frac{\vte-\overline\vt}\ep\Big)
\partial_t\Big(\alpha
R_\ep+\beta T_\ep\Big)}\ {\rm d}t
\end{equation}
$$
=
-\frac{\alpha}{\beta^2+\alpha\delta}\int_0^\tau\intO{\Big(\delta\frac{\vte-\overline\vt}\ep-\beta\frac{\vre-\overline\vr}\ep\Big)
\partial_t\Big(\delta T_\ep-\beta R_\ep\Big)}\ {\rm d} t
$$

Finally, the last line on the right-hand side of (\ref{*1}) reads
\begin{equation}\label{*4}
-\int_0^\tau\intO{\Big(\delta\frac{\vte-\Theta_\ep}\ep
-\beta\frac{\vr_\ep-r_\ep}\ep\Big)\vc U_\ep\cdot\Grad T_\ep}\ {\rm
d}t
\end{equation}
$$=
-\frac\beta{\beta^2+\alpha\delta}\int_0^\tau\intO{\Big(\delta\frac{\vte-\Theta_\ep}\ep
-\beta\frac{\vr_\ep-r_\ep}\ep\Big)\vc U_\ep\cdot\Grad\Big( \alpha
R_\ep +\beta T_\ep\Big)}\ {\rm d}t
$$
$$
-
\frac\alpha{\beta^2+\alpha\delta}\int_0^\tau\intO{\Big(\delta\frac{\vte-\Theta_\ep}\ep
-\beta\frac{\vr_\ep-r_\ep}\ep\Big)\vc U_\ep\cdot\Grad\Big( \delta
T_\ep -\beta R_\ep\Big)}\ {\rm d}t
$$
where the first term tends to zero due to the dispersion estimates (\ref{a2}).

Summing (\ref{*3}--\ref{*4}) we deduce the following result
$$
-\frac{\alpha}{\beta^2+\alpha\delta}
\int_0^\tau\intO{\Big(\delta\frac{\vte-\Theta_\ep}\ep
-\beta\frac{\vr_\ep-r_\ep}\ep\Big)\Big(\partial_t(\delta
T_\ep-\beta R_\ep) +\vc U_\ep\cdot\Grad (\delta T_\ep-\beta
R_\ep)\Big)}\ {\rm d}t
$$
$$
= \frac{\alpha}{\beta^2+\alpha\delta}
\int_0^\tau\intO{\Big(\delta\frac{\vte-\Theta_\ep}\ep
-\beta\frac{\vr_\ep-r_\ep}\ep\Big)\Delta\Phi_\ep}\ {\rm d}t+ \chi^1_{\ep}(\eta, \tau) = \chi_{\ep}(\eta, \tau)
,
$$
where we have used (\ref{r4a}), and, again, the dispersive
estimates (\ref{a2}). Resuming the calculations in this section, we can
rewrite inequality (\ref{r11+}) as follows
\bFormula{r12+} \left[ \mathcal{E}_\ep \left( \vre, \vte, \vue
\Big| r_\ep , \Theta_\ep , \vc{U}_\ep \right) \right]_{t = 0}^\tau
\leq \left[ \intO{ \Ov{\vr} \frac{1}{2}|\Grad \Phi_\ep |^2 }
\right]_{t = 0}^{t = \tau} \eF
\[
+ \frac12\frac\delta{\beta^2+\alpha\delta} \Big[\intO{|\alpha
R_\ep+\beta T_\ep|^2}\Big]_0^\tau+
\frac12\frac\alpha{\beta^2+\alpha\delta} \Big[\intO{|\delta
T_\ep-\beta R_\ep|^2}\Big]_0^\tau
\]
\[
+ \int_0^\tau c\left(1 + \chi^1_{\ep} (t, \eta) + \| \Grad T_\ep (t, \cdot) \|_{L^\infty(R^3;R^3)} \right) \mathcal{E}_\ep \left( \vre, \vte, \vue \Big| r_\ep , \Theta_\ep , \vc{U}_\ep \right) \ \dt   + \chi^2_\ep
(\tau, \eta),
\]
where
\[
\chi^i_{\ep}(\cdot, \eta) \to 0 \ \mbox{in} \ L^1(0,T) \ \mbox{as}\ \ep \to 0 \ \mbox{for any fixed}\ \eta > 0, \ i=1,2.
\]

Consequently, in accordance with the energy balances (\ref{a1}),
(\ref{a3}), { and dispersive estimates (\ref{a2})}, inequality
(\ref{r12+}) reduces to

\greybox{

\bFormula{r13}
\left[ \mathcal{E}_\ep \left( \vre, \vte, \vue
\Big| r_\ep , \Theta_\ep , \vc{U}_\ep \right) \right]_{t = 0}^\tau
\eF
\[
\leq \int_0^\tau c\left(1 + \chi^1_{\ep} (t, \eta) +
 \| \Grad T_\ep (t, \cdot) \|_{L^\infty(R^3;R^3)} \right) \mathcal{E}_\ep \left( \vre, \vte, \vue \Big| r_\ep , \Theta_\ep , \vc{U}_\ep \right) \ \dt   + \chi^2_\ep
(\tau, \eta),
\]
where
\[
\chi^i_{\ep}(\cdot, \eta) \to 0 \ \mbox{in} \ L^1(0,T) \ \mbox{as}\ \ep \to 0 \ \mbox{for any fixed}\ \eta > 0, \ i=1,2.
\]

}

\subsection{Conclusion}

Summarizing (\ref{b3}), (\ref{b4}) and (\ref{conv1}), we obtain
\[
\sqrt{\vre} \vue \to \Ov{\sqrt{\vr} \vu} \ \mbox{weakly in} \ L^\infty(0,T;L^2(R^3;R^3))
\]
\[
\frac{\vre - \Ov{\vr}}{\ep} = \left[ \frac{\vre - \Ov{\vr}}{\ep} \right]_{\rm ess} + \left[ \frac{\vre - \Ov{\vr}}{\ep} \right]_{\rm res},
\]
where
\[
\left[ \frac{\vre - \Ov{\vr}}{\ep} \right]_{\rm ess} \to R \ \mbox{weakly-(*) in} \ L^\infty(0,T; L^2(R^3)),
\]
while
\[
\left[ \frac{\vre - \Ov{\vr}}{\ep} \right]_{\rm res} \to 0 \ \mbox{in} \ L^\infty(0,T; L^{5/3}(R^3)).
\]

Similarly
\[
\left[ \frac{\vte - \Ov{\vt}}{\ep} \right]_{\rm ess} \to T \ \mbox{weakly-(*) in} \ L^\infty(0,T; L^2(R^3)),
\]
and
\[
\left[ \frac{\vte - \Ov{\vt}}{\ep} \right]_{\rm res} \to 0 \ \mbox{in} \ L^\infty(0,T; L^{q}(R^3)) \ \mbox{for any}\ 1 \leq q < 2.
\]

On the other hand,
\[
\Grad \Phi_\ep \to 0 \ \mbox{in} \ L^1(0,T; W^{k,p}(R^3;R^3)) \cap L^\infty_{\rm loc} ((0,T]; W^{k,p}(R^3;R^3)) \ \mbox{for any}\
2 < p \leq \infty,\ k=0,1,\dots,
\]
whereas
\[
R_\ep \to R_\eta, \ T_\ep \to T_\eta \ \mbox{in} \ L^1(0,T; W^{1,\infty} (R^3)) \cap L^\infty_{\rm loc} ((0,T]; W^{1, \infty}(R^3)),
\]
where, in view of the dispersive estimates (\ref{a2}),
\bFormula{eq1} \alpha R_\eta + \beta T_\eta = 0, \eF
and { due to (\ref{r4a}),}
\bFormula{eq2}
\partial_t (\delta T_\eta - \beta R_\eta) + \vc{v} \cdot \Grad (\delta T_\eta - \beta R_\eta) = 0
\eF
with the initial data
\[
R_{0, \eta} = \chi_\eta * [ \psi_\eta \vr^{(1)}_0 ] , \ T_{0, \eta} = \chi_\eta * [ \psi_\eta \vt^{(1)}_0 ].
\]

Now, applying Gronwall's lemma to (\ref{r13}) we obtain

\greybox{

\bFormula{gron}
\intO{ \left[ \frac{1}{2} \left| \sqrt{\vre} \vue - \sqrt{\vre} \Grad \Phi_\ep - \sqrt{\vre} \vc{v} \right|^2 (\tau, \cdot) \right] }
\eF
\[
+ \frac{1}{\ep^2} \intO{ \left[ H_{\Theta_\ep} (\vre, \vte) - \frac{\partial H_{\Theta_\ep} (r_\ep, \Theta_\ep) }{\partial \vr} \left(
\frac{\vre - \Ov{\vr}}{\ep} - R_\ep \right) - H_{\Theta_\ep} (r_\ep, \Theta_\ep )  \right] (\tau, \cdot) }
\]
\[
\leq \exp \left( \int_0^\tau c\left(1 + \chi^1_{\ep} (t, \eta) + \| \Grad T_\ep (t, \cdot) \|_{L^\infty(R^3;R^3)} \right) \dt \right)
\left[ \chi^2_\ep(\tau, \eta) + \frac{1}{2} \intO{ \vr_{0, \ep} \left| \vu_{0, \ep} - \Grad \Phi_{0,\ep} - \vc{v}_0 \right|^2 } \right]
\]
\[
+ c \exp \left( \int_0^\tau c\left(1 + \chi^1_{\ep} (t, \eta) + \| \Grad T_\ep (t, \cdot) \|_{L^\infty(R^3;R^3)} \right) \dt \right)
\left[ \left\| \vr^{(1)}_{0,\ep} - R_{0,\ep} \right\|^2_{L^2(R^3)} +  \left\| \vt^{(1)}_{0,\ep} - T_{0,\ep} \right\|^2_{L^2(R^3)} \right]
\]
for any $\tau \in [0,T]$.

}

Thus, letting $\ep \to 0$ in (\ref{gron}) and making use of the convergence relations established earlier in this section, we get
\bFormula{gron1}
\limsup_{\ep \to 0} \left( \int_K \left[ \frac{1}{2} \left| \sqrt{\vre} \vue  - \sqrt{\Ov{\vr}} \vc{v} \right|^2 (\tau, \cdot) \right] \dx \right.
\eF
\[
\left.
+
\frac{1}{\ep^2} \int_K \left[ H_{\Theta_\ep} (\vre, \vte) - \frac{\partial H_{\Theta_\ep} (r_\ep, \Theta_\ep) }{\partial \vr} \left(
\frac{\vre - \Ov{\vr}}{\ep} - R_\ep \right) - H_{\Theta_\ep} (r_\ep, \Theta_\ep )  \right] (\tau, \cdot)\ \dx \right)
\]
\[
\leq \exp \left( \int_0^\tau c\left(1 +  \| \Grad T_\eta (t, \cdot) \|_{L^\infty(R^3;R^3)} \right) \dt \right)
\left[ \frac{\Ov{\vr}}{2} \intO{  \Big| \Grad \Delta^{-1} [\Div \vu_{0}] - \Grad \left( \chi_\eta * (\psi_\eta \Delta^{-1} [\Div \vu_0]) \right) \Big|^2 } \right]
\]
\[
+ c \exp \left( \int_0^\tau c\left(1 +  \| \Grad T_\eta (t, \cdot) \|_{L^\infty(R^3;R^3)} \right) \dt \right)
\left[ \left\| \vr^{(1)}_{0} - \chi_\eta * (\psi_\eta \vr^{(1)}_0) \right\|^2_{L^2(R^3)} +  \left\| \vt^{(1)}_{0} -  \chi_\eta * (\psi_\eta \vt^{(1)}_0)\right\|^2_{L^2(R^3)} \right]
\]
for any $ \tau \in (0,T]$ and any compact $K \subset R^3$.

Finally, in accordance with (\ref{eq1}), (\ref{eq2}),
\bFormula{konec}
\partial_t T_\eta + \vc{v} \cdot \Grad T_\eta = 0,\ T_{\eta}(0, \cdot) =  \Ov{\vr} \frac{ \partial s(\Ov{\vr}, \Ov{\vt}) }{\partial \vt} \chi_\eta * [ \psi_\eta \vt^{(1)}_0 ] -
\frac{1}{\Ov{\vr}} \frac{ \partial p(\Ov{\vr}, \Ov{\vt}) }{\partial \vt} \chi_\eta * [ \psi_\eta \vr^{(1)}_0 ]  ;
\eF
whence, by virtue of hypothesis (\ref{HYP}),
\[
\| \Grad T_\eta \|_{L^\infty(R^3;R^3)} \ \mbox{is bounded in}\ L^\infty(0,T) \ \mbox{uniformly for}\ \eta \to 0.
\]
Consequently, making use of the estimate
\[
\frac{1}{\ep^2} \int_K \left[ H_{\Theta_\ep} (\vre, \vte) - \frac{\partial H_{\Theta_\ep} (r_\ep, \Theta_\ep) }{\partial \vr} \left(
\frac{\vre - \Ov{\vr}}{\ep} - R_\ep \right) - H_{\Theta_\ep} (r_\ep, \Theta_\ep )  \right] (\tau, \cdot)\ \dx
\]
\[
\geq c \left( \left\| \left[ \frac{\vre - \Ov{\vr}}{\ep} - R_\ep \right]_{\rm ess} \right\|^2_{L^2(K)} +
\left\| \left[ \frac{\vte - \Ov{\vt}}{\ep} - T_\ep \right]_{\rm ess} \right\|^2_{L^2(K)} \right)
\]
we may let $\eta \to 0$ in (\ref{gron1}) to obtain the desired conclusion (\ref{m4}), (\ref{m5}). Thus, passing to the limit $\eta \to 0$ in
(\ref{konec}) and letting $\ep \to 0$ in the momentum equation (\ref{v2}) for solenoidal test functions $\varphi$ completes the proof of Theorem \ref{Tm1}.

\section{Concluding remarks}

Similar results can be obtained on a general (unbounded) domain $\Omega \subset R^3$ as soon as the following conditions hold:

\begin{itemize}
\item
the velocity field $\vc{u}$ satisfies the \emph{complete slip} conditions
\[
\vc{u} \cdot \vc{n}|_{\partial \Omega} = 0,\ ( \tn{S} (\vt, \Grad \vu) \vc{n} ) \times \vc{n} |_{\partial \Omega} = 0,
\]
or \emph{Navier's boundary conditions}
\[
\vc{u} \cdot \vc{n}|_{\partial \Omega} = 0,\ [\tn{S} (\vt, \Grad \vu) \vc{n}]_{\rm tan}  + \beta \vc{u} |_{\partial \Omega} = 0,
\]
where $\beta \geq 0$ is a ``friction'' coefficient;

\item the target Euler system (\ref{l1}), (\ref{l2}) possesses a regular solution on $[0, T_{\rm max})$;

\item the acoustic equation (\ref{r3}), (\ref{r4}) admits the dispersive estimates (\ref{a2});

\item the generalized Korn inequality holds: For any $M >0$ there exists $c(M) > 0$ such that
\[
\| \vc{w} \|^2_{W^{1,2}(\Omega;R^3)} \leq c(M) \left( \left\| \Grad \vc{w} + \Grad^t \vc{w} - \frac{2}{3} \Div \vc{w} \tn{I} \right\|^2_{L^2(\Omega; R^{3 \times 3})} + \int_{ \Omega \setminus V } |\vc{w}|^2 \ \dx \right) ,\ \vc{w} \in W^{1,2}(\Omega;R^3),
\]
for any measurable set $V \subset \Omega$, $|V| \leq M$.

\end{itemize}

These conditions are satisfied, for example, if $\Omega \subset R^3$ is an \emph{exterior} domain with Lipschitz boundary, see Alazard \cite{AL2}, Isozaki \cite{Isoz}, and \cite[Appendix]{FeNo6}.


\def\ocirc#1{\ifmmode\setbox0=\hbox{$#1$}\dimen0=\ht0 \advance\dimen0
  by1pt\rlap{\hbox to\wd0{\hss\raise\dimen0
  \hbox{\hskip.2em$\scriptscriptstyle\circ$}\hss}}#1\else {\accent"17 #1}\fi}

\end{document}